\documentclass[11pt]{article}
\usepackage{latexsym}
\usepackage[all]{xy}
\usepackage{amsmath}
\usepackage{amssymb}
\usepackage{amsfonts}
\usepackage[applemac]{inputenc}

\newtheorem{thm}{Theorem}[section]
\newtheorem{prop}[thm]{Proposition}
\newtheorem{lem}[thm]{Lemma}

\newtheorem{df}[thm]{Definition}
\newtheorem{cor}[thm]{Corollary}

\newtheorem{rmk}[thm]{Remark}

\newcommand{\cdga}{\mathsf{cdga}_{k}}
\newcommand{\dglie}{\mathsf{dgLieAlg}_{k}}

\newcommand{\dgliea}{\mathsf{dgLieAlgd}}

\begin{document}

\title{\textbf{A model structure on relative dg-Lie algebroids}}
\bigskip
\bigskip
\bigskip
\bigskip

\author{
\bigskip \\ \bf{Gabriele Vezzosi}\\
\small{Institut de Mathématique de Jussieu} \\ 
\small{Paris - France}
\\}

\date{June 2014}

\maketitle

\begin{abstract}
In this paper, for the future purposes of relative formal derived deformation theory, we prove the existence of a model structure on the category of dg-Lie algebroids over a cochain differential non-positively  graded commutative algebra over a commutative base ring of characteristic zero.  
\end{abstract}


\bigskip
\bigskip
\bigskip

\bigskip

\section{Introduction} This paper provides a model structure on the category of dg-Lie algebroid over a non positivey graded cdga over a base ring $k$ of characteristic zero. This result might be useful in various contexts but it was  originally motivated from the need of such a structure in \emph{formal deformation theory} relative to an arbitrary cdga $A$ (\cite{benj, tvv}), extending the theory relative to $A=k$ developed in \cite{lufor}. In this more general situation, the role of dg-Lie algebras is played by dg-Lie algebroids, hence the connection with the result of this paper.\\

Again for the sake of motivating the result in this paper, let me sketch briefly what this relative deformation theory is about. Details will appear in \cite{tvv}. First of all, the relative formal deformation theory I am alluding at involves the study of formal deformation of \emph{tame} dg-categories.  The adjective \emph{tame} here refers to a model structure on the category of dg-categories over $A$ (i.e categories enriched in $A$-dg-modules). This model structure is induced by the so-called \emph{tame homological algebra}, or more precisely by the tame model structure on $A$-dg-modules, which is defined as follows.  Let $\mathcal{P}_{A}$ be the class of all 
$A$-dg-modules which are graded projectives as graded $A$-modules (i.e. retracts of direct sums of shifts of $A$). Then we define a map $f : E \longrightarrow F$ of $A$-dg-modules to be a  \emph{tame quasi-isomorphism}
if for all $P \in \mathcal{P}$ the induced morphism of complexes of abelian groups
$\underline{\mathsf{Hom}}(P,E) \longrightarrow \underline{\mathsf{Hom}}(P,F)$
is a quasi-isomorphism. The tame quasi-isomorphisms together with the projective fibrations (i.e. epimorphisms) define the tame model structure on the category of $A$-dg-modules. This model structure is stable and combinatorial. Moreover one can show that the usual projective model structure is a right Bousfield localization of the tame one (\cite{tvv}). In particular, the tame equivalence are projective equivalences i.e. quasi-isomorphisms, and  moreover homotopy equivalences are tame equivalences.
Therefore, as explained in Remark \ref{tame}, the result of this paper also extends to give a \emph{tame} model structure on dg-Lie algebroids over $A$. This will prove useful in \cite{tvv}.\\

\noindent \textbf{Notations.} 

\begin{itemize}
\item  $k$ will denote the base field and will be assumed to be of characteristic $0$.
\item $\mathsf{C}(k)$ will denote the model category of unbounded complexes of $k$-modules with surjections as fibrations, and quasi-isomorphisms s equivalences. It is a symmetric monoidal model category with the usual tensor product $\otimes_{k}$ of complexes over $k$.
\item $\cdga$ denotes the category of differential \emph{nonpositively} graded algebras over $k$, with differential increasing the degree by 1. 
We will always consider $\cdga$ endowed with the usual model structure for which fibrations are surjections in negative degrees, and equivalences are quasi-isomorphisms (see \cite[\S 2.2.1]{hagII}).
\item  A dg-Lie algebra over $k$ is a Lie algebra object in the symmetric monoidal category $(\mathsf{C}(k), \otimes_{k})$, i.e. an object $\mathfrak{g}$ in $\mathsf{C}(k)$ endowed with a morphism of complexes (called the bracket) $[-,-]:\mathfrak{g}\otimes \mathfrak{g} \rightarrow \mathfrak{g}$ such that 
\begin{enumerate}
\item $[x,y]= (-1)^{|x||y| +1} [y,x]$
\item $[x, [y,z]] + (-1)^{(|z||x| + |x||y|)}[y,[z,x]] + (-1)^{(|z||x| + |z||y|)}[z,[y,x]]=0$
\item $d[x,y]= [dx,y] + (-1)^{|x|}[x,dy]$
\end{enumerate}  
A morphism of dg-Lie algebras over $k$ is a morphism of complexes commuting with the brackets.
\item We'll write $\dglie$ for the category of unbounded dg-Lie algebras over $k$. By \cite[Thm 4.1.1]{hin} (with no need of its erratum, since $k$ has characteristic zero), applied to the $k$-dg-operad $\mathcal{O}=\mathsf{Lie}$, $\dglie$ is a (left proper, combinatorial) 
model category of unbounded dg-Lie algebras over $k$, where fibrations and equivalences are defined via the forgetful functor $\dglie \rightarrow \mathsf{C}(k)$. Moreover, $\dglie$ is a simplicial model category via $\underline{\mathrm{Hom}}(\mathfrak{g}, \mathfrak{h})_{n}:= \mathrm{Hom}_{\dglie}(\mathfrak{g}, \Omega^{\bullet}(\Delta_{n})\otimes_{k} \mathfrak{h})$.
\end{itemize}
\section{Relative dg-Lie algebroids}
We start by defining what is a dg-Lie algebroid over a cdga $A$.

\begin{df} \begin{itemize}
\item Let $A\in \cdga$. The category $\dgliea _{A}$ of \emph{dg-Lie algebroids over} $A$ is the category whose
\begin{itemize}
\item \emph{objects} are the pairs $(L, \alpha: L\rightarrow T_{A})$ where $L$ is a $A$-dg-module and a $k$-dg-Lie algebra, and $\alpha$ (called the anchor map) is a morphism of $A$-dg-modules and of $k$-dg-Lie algebras, such that for any homogeneous $\ell_1, \ell_2 \in L$ and any homogeneous  $a\in A$, the following \emph{graded Leibniz rule} holds $$[\ell_1, a\ell_2]= (-1)^{|a||\ell_{1}|}a[\ell_1,\ell_2]+\alpha(\ell_1)(a)\ell_2$$
\item \emph{morphisms} $(L', \alpha')\rightarrow (L,\alpha)$  are morphisms $\psi: L'\rightarrow L$ of $A$-dg-modules and of $k$-dg-Lie algebras, commuting with anchor maps. A morphism in $\dgliea_{A}$ is an \emph{equivalence} if it is a quasi-isomorphism of the underlying $A$-dg modules.
\end{itemize}
\end{itemize}
\end{df}

\noindent \textbf{Examples.} When $A=k$ we get the category of  dg-Lie $k$-algebras (since $T_A$ is trivial in this case, so the anchor map is uniquely determined).  $(T_{A}, \textrm{id}: T_{A}\rightarrow T_{A})$ is the final object in $\dgliea _{A}$, and $\mathrm{Hom}_{\dgliea _{A}}((T_{A}, \textrm{id}: T_{A}\rightarrow T_{A}), (L, \alpha: L\rightarrow T_{A}))$ is the set of integrable connections on the dg-Lie algebroid $(L, \alpha: L\rightarrow T_{A})$ (see \cite{bb} for the case where $A$ is discrete). 

\begin{rmk} \emph{From the point of view of derived algebraic geometry, the category $\dgliea _{A}$ is reasonable only for \emph{cofibrant} $A$, because in this case $T_A$ is indeed the tangent complex of $A$.}
\end{rmk}

In order to prove our main result, we will make use of the following \emph{transfer criterion} originally due to Quillen (\cite[II.4]{qui}):

\begin{prop} \label{transf} \emph{\cite[2.5, 2.6]{bm}} Let $\mathcal{D}$ be a cofibrantly generated model category and let $F : \mathcal{D}  \to \mathcal{E} $ be a left adjoint with right adjoint $G$. Assume that $\mathcal{E}$ has small colimits and finite limits. Define a map $f$ in $\mathcal{E}$ to be a weak equivalence (resp. fibration) iff $G(f)$ is a weak equivalence (resp. fibration). Then this defines a cofibrantly generated model structure on $\mathcal{E}$ provided
\begin{itemize}
\item $F$ preserves small objects
\item $\mathcal{E}$ has a fibrant replacement functor
\item $\mathcal{E}$ functorial path-objects for fibrant objects.
\end{itemize}
\end{prop}

\begin{thm} Let $A\in \cdga$ be cofibrant. Then the category $\dgliea _{A}$ of dg-Lie algebroids over $A$, endowed with equivalences and fibrations defined on the underlying $A$-dg-modules, is a cofibrantly generated model category.
\end{thm}

\noindent \textbf{Proof.} We will transfer the model structure from $A-\mathsf{dgMod}/T_{A}$ to $\dgliea_{A}$ via the adjunction $$\xymatrix{A-\mathsf{dgMod}/T_{A} \ar@<+ .7ex>[r]^-{\textrm{Free}} & \dgliea_{A} \ar@<+ .7ex>[l]^-{\textrm{Forget}}}$$ where the Free functor, i.e. the left adjoint, is defined in \cite{ka}. In order to do this, since $A-\mathsf{dgMod}/T_{A}$ is a cofibrantly generated model category and $\dgliea_{A}$ has small colimits and finite limits, we may use the transfer criterion of Prop. \ref{transf}, i.e. \\

\noindent \emph{Transfer criterion}. Define a map $f$ in $\dgliea_{A}$ to be a weak equivalence (resp. a fibration) if $\textrm{Forget}(f)$ is a weak equivalence (resp. a fibration). This defines a cofibrantly generated model structure on $\dgliea_{A}$, provided\\

\noindent (i) the functor $\textrm{Forget}$ preserves small objects;\\
\noindent (ii) $\dgliea_{A}$ has a fibrant replacement functor, and functorial path-objects for fibrant objects.\\

\noindent Now, condition (i) follows easily from the fact that the forgetful functor commutes with filtered colimits. For the first part of condition (ii) we may take the identity functor as a fibrant replacement functor, since all objects in $\dgliea_{A}$ are fibrant. To prove the second half of condition (ii) - i.e. the existence of functorial path-objects - we consider the following construction. Recall first that if $L_1 \rightarrow T_{A}$ and $L_2 \rightarrow T_A$ are dg-Lie algebroids over $A$, then their product in $\dgliea_{A}$ is $L_1 \times_{T_{A}}L_2 \rightarrow T_{A}$, with the obvious induced $k$-dg-Lie bracket. \\ Let now $(\alpha: L \rightarrow T_{A}) \in \dgliea_{A}$, and consider the $A$-dg-module $L[t,dt] := L\otimes_{k}k[t,dt]$, where $k[t,dt]:= k[t] \oplus k[t]dt$ is identified with the de Rham commutative dg-algebra of $\mathrm{Spec}(k[x_0,x_1]/(x_0 + x_1 -1))$ with $k[t]$ sitting in degree $0$, and $k[t]dt$ in degree $1$, and the differential is given by $$d(\ell_{1}f(t) + \ell_{2}g(t)dt):=d_{L}(\ell_1)f(t) + ((-1)^{|\ell_{1}|}\ell_1 f'(t)+ d_{L}(\ell_{2})g(t))dt.$$ We endow $L[t,dt]$ with the following $k$-dg-Lie bracket 
$$[\ell_{1}f_1(t), \ell_{2}f_{2}(t)]:= [\ell_1, \ell_2]_{L}f_1(t)f_2(t), \qquad [\ell_1 f_1(t), \ell_{2}g_{2}(t)dt]:= [\ell_1, \ell_2]_{L}f_1(t)g_2(t)dt.$$
Note that for any $s\in k$, we have an evaluation map $$\mathrm{ev}_{s}: L[t,dt] \longrightarrow L \, , \,\,\,\, \ell_1f(t) + \ell_2 g(t)dt \longmapsto f(s)\ell_{1}$$ that is a morphism of $A$-dg-modules and of $k$-dg-Lie algebras, and a quasi-isomorphism. \\
Consider now the pull-back diagram in $A-\mathsf{dgMod}$ (defining the underlying $A$-dg module of what will be the path-object $\textrm{Path}(L \rightarrow T_{A})$ in $\dgliea_{A}$)
$$\xymatrix{\textrm{Path}(L \rightarrow T_{A}) \ar[r] \ar[d] & L[t,dt] \ar[d] \\ T_{A} \ar[r] & T_{A}[t,dt] },$$ where the right vertical map is given by $$ L[t,dt]  \longrightarrow T_{A}[t,dt] \, , \,\,\ \sum_{i} \xi_{i}f_i(t) + \eta_{i}g_i (t) dt \longmapsto \sum_{i} \alpha(\xi_{i})f_i(t) + \alpha(\eta_{i})g_i (t) dt ,$$ while the map $T_{A} \rightarrow T_{A}[t,dt]$ is the obvious inclusion $\xi \mapsto \xi \cdot 1 + 0 \cdot dt$. In other words, $\textrm{Path}(\xymatrix{L \ar[r]^-{\alpha} & T_{A}})$ is the sub-$A$-dg module of $L[t,dt]$ consisting of elements $\sum_{i} \xi_{i}f_i(t) + \eta_{i}g_i (t) dt$ such that $$\sum_{i} \alpha(\xi_{i})f_i(t) \in  T_{A} \hookrightarrow T_{A}[t] \, , \,\,\,\,\,  \sum_{i}  \alpha(\eta_{i})g_i (t) dt=0 \, \, \textrm{in} \,\, T_{A} \hookrightarrow T_{A}[t]dt. $$ It is straightforward that the composition 
$$\xymatrix{\textrm{Path}(L \rightarrow T_{A})\ar[r] & L[t,dt] \ar[rr]^-{(\textrm{ev}_{0}, \textrm{ev}_{1})}  & & L\times L}$$ factors through the inclusion $L\times_{T_{A}} L \hookrightarrow L \times L$, and that the resulting diagram $$\xymatrix{\textrm{Path}(L \rightarrow T_{A}) \ar[rr]^-{p} \ar[rd] &  & L\times_{T_{A}} L \ar[ld] \\ & T_{A} & }$$ commutes. Moreover, $\textrm{Path}(\xymatrix{L \ar[r]^-{\alpha} & T_{A}})$ is closed under the $k$-dg-Lie bracket in $L[t,dt]$, and the canonical map $\textrm{Path}(\xymatrix{L \ar[r]^-{\alpha} & T_{A}}) \rightarrow T_{A}$ is a morphism of $k$-dg-Lie algebras. The fact that  $(\textrm{Path}(\xymatrix{L \ar[r]^-{\alpha} & T_{A}}) \rightarrow T_{A})$ is actually a dg-Lie algebroid over $A$,  i.e. satisfies also the graded Leibniz rule, follows from a direct computation. Namely, the lhs of the Leibniz identity is
$$[ \sum_{i} \xi_{i}f_i(t) + \eta_{i}g_i (t) dt, a(\sum_{j} \overline{\xi}_{j}\overline{f}_j(t) + \overline{\eta}_{j}\overline{g}_j (t) dt) ]=$$ $$= \sum_{i, j} (-1)^{|a||\xi_{i}|} a[\xi_{i}, \overline{\xi}_{j}] f_i(t)\overline{f}_{j}(t) + (-1)^{|a||\xi_{i}|} a[\xi_{i}, \overline{\eta}_{j}] f_i(t)\overline{g}_{j}(t)dt + (-1)^{|a||\eta_{i}| + (|a| + |\overline{\xi}_{j}|)} a[\eta_{i}, \overline{\xi}_{j}] \overline{f}_{j}(t)g_{i}(t)dt +$$ $$+  \alpha(\xi_{i})(a)\overline{\xi}_{j} f_{i}(t)\overline{f}_{j}(t) + \alpha(\xi_{i})(a)\overline{\eta}_{j} f_{i}(t)\overline{g}_{j}(t)dt + (-1)^{|a| + |\overline{\xi}_{j}|} \alpha(\eta_{i})(a)\overline{\xi}_{j} \overline{f}_{j}(t)g_{i}(t)dt ;$$ we observe that, since $(\sum_{i} \xi_{i}f_i(t) + \eta_{i}g_i (t) dt) \in \textrm{Path}(\xymatrix{L \ar[r]^-{\alpha} & T_{A}})$, the last term $\alpha(\eta_{i})(a)\overline{\xi}_{j} \overline{f}_{j}(t)g_{i}(t)dt$ vanishes, and moreover $\sum_{i} \alpha(\xi_{i}) f_{i}(t) + \alpha(\eta_{i})g_{i}(t)dt = \partial $, for some $\partial \in T_{A}$, so that 
$$[ \sum_{i} \xi_{i}f_i(t) + \eta_{i}g_i (t) dt, a(\sum_{j} \overline{\xi}_{j}\overline{f}_j(t) + \overline{\eta}_{j}\overline{g}_j (t) dt) ]=$$ $$= \sum_{i, j} ((-1)^{|a||\xi_{i}|} a[\xi_{i}, \overline{\xi}_{j}] f_i(t)\overline{f}_{j}(t) + (-1)^{|a||\xi_{i}|} a[\xi_{i}, \overline{\eta}_{j}] f_i(t)\overline{g}_{j}(t)dt + (-1)^{|a||\eta_{i}| + (|a| + |\overline{\xi}_{j}|)} a[\eta_{i}, \overline{\xi}_{j}] \overline{f}_{j}(t)g_{i}(t)dt ) +$$ $$+ \sum_{j} \partial(a) \overline{\xi}_{j}\overline{f}_{j}(t) + \partial(a) \overline{\eta}_{j}\overline{g}_{j}(t)dt.$$

The rhs of the Leibniz identity reads
$$\sum_{i,j}[\xi_{i}f_{i}(t) , a(\sum_{j} \overline{\xi}_{j}\overline{f}_j(t) + \overline{\eta}_{j}\overline{g}_j (t) dt)] + [\eta_{i}g_{i}(t)dt , a(\sum_{j} \overline{\xi}_{j}\overline{f}_j(t) + \overline{\eta}_{j}\overline{g}_j (t) dt))] =$$ $$= \sum_{i,j} ((-1)^{|a||\xi_{i}|} a[\xi_{i}f_{i}(t), \overline{\xi}_{j}\overline{f}_j(t) + \overline{\eta}_{j}\overline{g}_j (t) dt) ] + (-1)^{|a|(|\eta_{i} +1|} a[\eta_{i}g_{i}(t)dt, \overline{\xi}_{j}\overline{f}_j(t) + \overline{\eta}_{j}\overline{g}_j (t) dt) ]) +$$ $$ + \partial(a)(\sum_{j} \overline{\xi}_{j}\overline{f}_j(t) + \overline{\eta}_{j}\overline{g}_j (t) dt) =$$ $$=  \sum_{i, j} ((-1)^{|a||\xi_{i}|} a[\xi_{i}, \overline{\xi}_{j}] f_i(t)\overline{f}_{j}(t) + (-1)^{|a||\xi_{i}|} a[\xi_{i}, \overline{\eta}_{j}] f_i(t)\overline{g}_{j}(t)dt + (-1)^{|a||\eta_{i}| + (|a| + |\overline{\xi}_{j}|)} a[\eta_{i}, \overline{\xi}_{j}] \overline{f}_{j}(t)g_{i}(t)dt ) +$$ $$ + \sum_{j} \partial(a) \overline{\xi}_{j}\overline{f}_{j}(t) + \partial(a) \overline{\eta}_{j}\overline{g}_{j}(t)dt,$$ and we conclude that  $(\textrm{Path}(\xymatrix{L \ar[r]^-{\alpha} & T_{A}}) \rightarrow T_{A})$ is indeed a dg-Lie algebroid over $A$.\\

Finally, the diagonal map of dg-Lie algebroids over $A$ $$L\longrightarrow L\times_{T_{A}}L$$ factors as $$\xymatrix{L \ar[r]^-{u} & \textrm{Path}(L \rightarrow T_{A}) \ar[r]^-{p} & L\times_{T_{A}}L},$$ where $u(\xi)= \xi \cdot 1 + 0\cdot dt$. Now, $p$ is surjective, since for any $(\xi,\eta) \in L\times_{T_{A}}L$, we have that $\xi t + \eta (1-t) \in \textrm{Path}(L \rightarrow T_{A})$ and $p(\xi t + \eta (1-t))=(\xi,\eta)$; hence $p$ is a fibration. The following Lemma shows the existence of path-objects in  $\dgliea_{A}$.

\begin{lem}\label{uqiso} The map $u:L \longrightarrow \textrm{Path}(L \rightarrow T_{A})$ is a weak equivalence (i.e. a quasi-isomorphism).
\end{lem}

\noindent \textsf{Proof of lemma.} First of all, we notice that the evaluation-at-$0$ map $$\textrm{ev}_{0}: \textrm{Path}(L \rightarrow T_{A}) \longrightarrow L  \, , \,\,\,  \sum_{i} \xi_{i}f_i(t) + \eta_{i}g_i (t) dt \longmapsto \sum_{i} \xi_{i}f_i(0)$$ is a left inverse to $u$. So it will be enough to produce a homotopy between $u\circ \textrm{ev}_{0}$ and $\textrm{Id}_{\textrm{Path}(L \rightarrow T_{A})}$. Let us define a family of $k$-linear maps, indexed by $p\in \mathbb{Z}$, $$h^{p}:\textrm{Path}(L \rightarrow T_{A})^{p} \longrightarrow \textrm{Path}(L \rightarrow T_{A})^{p-1} \, , \,\,\,\,  \xi f(t) + \eta  g(t) dt \longmapsto (-1)^{p} \eta \int_{0}^{t}g(x)dx. $$
We pause a moment to explain why this maps are well defined. We have used that if $\theta: =\sum_{i} \xi_{i}f_i(t) + \eta_{i}g_i (t) dt $ is homogeneous of degree $p$ (so that $|\xi_{i}|= p$ and $|\eta_{i}|= p-1$, for all $i$), then $h^{p}(\theta)$, which a priori belongs to $L[t,dt]^{p-1}$, actually belongs to $\textrm{Path}(L \rightarrow T_{A})^{p-1}$. This can be seen as follows. Write $\sum_{i} \eta_{i}g_i (t) dt$ as $\sum_{s\geq 0}\tilde{\eta}_{s}t^{s} dt$, where each $\tilde{\eta}_{s}$ is a $k$-linear combination of all the $\eta_{i}$'s. Then the condition that $\theta \in \textrm{Path}(L \rightarrow T_{A})^{p}$ reads $\sum_{s\geq 0} \alpha(\tilde{\eta}_{s})t^{s} dt= 0$, i.e. that $\alpha(\tilde{\eta}_{s})=0$ for all $s$. Therefore $$h^{p}(\theta)= \sum_{s\geq 0}\tilde{\eta}_{s}t^{s+1}/(s+1) dt$$ belongs to $\textrm{Path}(L \rightarrow T_{A})^{p-1}$, as claimed.\\
Now, it is a straightforward verification that the family of maps $\{h^{p}\}$ yields a homotopy between $u\circ \textrm{ev}_{0}$ and $\textrm{Id}_{\textrm{Path}(L \rightarrow T_{A})}$, i.e. that $$(h^{p+1}\circ d +d \circ h^{p})(\xi f(t) + \eta  g(t) dt)= -\xi (f(t)-f(0)) -\eta g(t)dt = (u\circ \textrm{ev}_{0} - \textrm{Id}_{\textrm{Path}(L \rightarrow T_{A})}) (\xi f(t) + \eta  g(t) dt).$$

If $\psi: (L, \alpha)\rightarrow (L',\alpha')$  is a morphism in  $\dgliea_{A}$, then by definition of the $A$-dg modules $\textrm{Path}(\alpha: L \rightarrow T_{A})$ and $\textrm{Path}(\alpha': L' \rightarrow T_{A})$, there is an induced map of $A$-dg modules $$\textrm{Path}(\psi): \textrm{Path}(L \rightarrow T_{A}) \longrightarrow \textrm{Path}(L' \rightarrow T_{A})$$ which obviously commutes with the anchor maps of the corresponding dg-Lie algebroid structures. And it is easy to verify that:
\begin{itemize}
\item $\textrm{Path}(\psi)$ is also a map of $k$-dg Lie algebras, so it is a morphism in $\dgliea_{A}$ 
\item if $\rho: (L', \alpha')\rightarrow (L'',\alpha'')$  is a morphism in  $\dgliea_{A}$, then $\textrm{Path}(\rho \circ \psi) = \textrm{Path}(\rho) \circ \textrm{Path}(\psi)$, and $\textrm{Path}(\mathrm{id}_{(L, \alpha)})= \mathrm{id}_{\textrm{Path}(\alpha: L \rightarrow T_A)}$
\item the following diagram (with obvious notations from Lemma \ref{uqiso}) commutes $$\xymatrix{ L \ar[rr]^-{\psi} \ar[d]_-{u}  & & L' \ar[d] ^-{u'} \\  \textrm{Path}(L \rightarrow T_{A}) \ar[rr]_-{\textrm{Path}(\psi)} & &  \textrm{Path}(L \rightarrow T_{A})}$$
\end{itemize}

In other words, the assignment $(\alpha: L \rightarrow T_{A}) \longmapsto \textrm{Path}(\alpha: L \rightarrow T_{A})$ is indeed a \emph{functorial} path object in $\dgliea_{A}$. Therefore, by Prop. \ref{transf}, we deduce the existence of the transferred model structure on $\dgliea_{A}$.

\hfill $\Box \,\,\,\, \Box$

As a consequence we obtain the following fact, first proved by V. Hinich (\cite[Thm. 4.1.1]{hin})

\begin{cor} For any discrete $\mathbb{Q}$-algebra $k$, the category $\dglie$ is a model category with weak equivalences and fibrations detected by the forgetful functor to $\mathsf{C}(k)$.
\end{cor}

\begin{rmk}\label{tame} \emph{Exactly the same argument proves that we may transfer the \textit{tame} model structure (as defined in the Introduction) from $A-\mathsf{dgMod}/T_{A}$ to $\dgliea_{A}$ via the same adjunction. In fact, fibrations remain the same, and Lemma \ref{uqiso} is still true since a homotopy equivalence is a tame equivalence.}
\end{rmk}

\end{document}